\documentclass[12pt,reqno]{amsart}
\usepackage{amssymb,latexsym,amsmath,amsthm,amscd,epsfig}
\usepackage[latin1]{inputenc}

\def\cG{\mathcal{G}}

\def\bpm{\begin{pmatrix}}
\def\epm{\end{pmatrix}}

\newcommand{\rf}[1]{{\eqref{#1}}}

\newcommand{\fiproof}{{\hspace*{\fill} $\square$ \vspace{2pt}}}
\newcommand{\D}{{\mathbb D}}

\newtheorem{thm}{Theorem}[section]

\newtheorem{lemma}[thm]{Lemma}

\newtheorem{theorem}[thm]{Theorem}

\theoremstyle{remark}

\numberwithin{equation}{section}

\def\C{\mathbb C}
\def\D{\mathbb D}

\def\cB{\mathcal B}
\def\cG{\mathcal G}

\def\diam{\text{diam}}

\def\diam{\operatorname{diam}}

\newcommand{\ve}{{\varepsilon}}

\title{Quasiconformal distortion of Hausdorff measures}

\author{Xavier Tolsa}
\address{Instituci\'{o} Catalana de Recerca i Estudis Avan\c{c}ats (ICREA) and Departament de Matem\`{a}tiques,  Universitat Aut\`{o}noma de Barcelona, 08193 Bellaterra (Barcelona), Catalonia}
\email{{\tt xtolsa@math.uab.cat}}
\urladdr{http://mat.uab.es/~xtolsa}

\thanks{Partially supported by grant MTM2007-62817 (Spain) and
2009-SGR-420 (Catalonia).}

\begin{document}

\maketitle

\begin{abstract}
In this paper we prove that if $\phi:\C\to\C$ is a $K$-quasiconformal map, $0<t<2$,
and $E\subset \C$ is a compact set contained in a ball $B$, then
$$\frac{H^t(E)}{\diam(B)^t} \leq C(K) \,\left(\frac{H^{t'}(\phi(E))}{\diam(\phi(B))^{t'}}\right)^{\frac{t}{t'K}},$$ 
where $t' = \frac{2t}{2K-Kt+t}$ and $H^s$ stands for the $s$-Hausdorff measure. 
In particular, this implies that $\phi$ transforms sets of finite $t'$-Hausdorff measure
into sets of finite $t$-Hausdorff measure.
\end{abstract}

\thispagestyle{empty}


\section{Introduction}

A homeomorphism $\phi:\Omega\to \Omega'$ between planar domains is called $K$-quasiconformal
if it preserves orientation, it belongs to the Sobolev space $W^{1,2}_{\rm loc}(\Omega)$, and it satisfies
$$\max_\alpha|\partial_\alpha\phi| \leq K\,\min_\alpha|\partial_\alpha\phi|\quad \mbox{ a.e. in $\Omega$.}$$  When $K=1$, the class of quasiconformal maps coincides with the
one of conformal maps.

 In his celebrated paper on distortion of area \cite{Astala-areadistortion}, 
 Astala proved that if $E\subset\C$ is compact, $\dim_H(E)$ denotes its Hausdorff dimension, and $\phi$ is $K$-quasiconformal, then
$$\frac1K \biggl(\frac1{\dim_H(E)} - \frac12 \biggr) \leq
\frac1{\dim_H(\phi(E))} - \frac12 \leq K \biggl(\frac1{\dim_H(E)} - \frac12 \biggr).$$

In a recent remarkable work, Lacey, Sawyer and Uriarte-Tuero \cite{Lacey-Sawyer-Uriarte} 
showed that if $E$ has positive
$t$-dimensional Hausdorff measure, then $\phi(E)$ has also positive $t'$-Hausdorff
measure, with $t'=\frac{2t}{2K-Kt+t}$ (the relationship between $t$ and $t'$ corresponds to the endpoint case 
in the rightmost inequality above).
Indeed, they proved that if $H^s_ \infty$ stands for the $s$-dimensional Hausdorff 
content and $E$ is contained in a ball $B$, then
\begin{equation}\label{eqcontent}
\frac{H_\infty^t(E)}{\diam(B)^t} \leq C(K) \,\left(\frac{H_\infty^{t'}(\phi(E))}{\diam(\phi(B))^{t'}}\right)^{\frac{t}{t'K}}.
\end{equation}
The case $t=1$ of \rf{eqcontent} had been proved previously in \cite{ACMOU}.
On the other hand, among other results,  in \cite{ACTUV} it has been shown that if
$H^{t'}(E)$ is $\sigma$-finite, then $H^{t}(\phi(E))$ is also
$\sigma$-finite (here $H^s$ stands for the $s$-Hausdorff measure). 

In this paper we will show that if $H^{t'}(E)$ is finite, then $H^{t}(\phi(E))$ is also
finite. Moreover, we will prove that the estimate \rf{eqcontent} also holds if we
replace Hausdorff contents by Hausdorff measures. The precise result is the following.

\begin{theorem} \label{teopri}
 Let $0<t<2$ and denote 
 $t' = \frac{2t}{2K-Kt+t}$. Let $\phi:\C\to\C$ be \hbox{$K$-quasiconformal}. For any compact set $E\subset\C$, we have
 \begin{equation}\label{eqhaus}
\frac{H^t(E)}{\diam(B)^t} \leq C(K) \,\left(\frac{H^{t'}(\phi(E))}{\diam(\phi(B))^{t'}}\right)^{\frac{t}{t'K}}.
\end{equation}
\end{theorem}

Notice that if we apply the preceding result to the map $\phi^{-1}$ (which is also $K$-quasiconformal) and to the set 
$\phi(E)$, we get
$$H^{t'}(E)<\infty\quad \Rightarrow \quad H^{t}(\phi(E))<\infty.$$
When $t=2$, the theorem also holds: it 
coincides with Astala's theorem on the distortion of area \cite{Astala-areadistortion}
(see also \cite[p.325]{AIM}).

To prove Theorem \ref{teopri} we will use a suitable 
version of the classical Frostman lemma and a result 
from \cite{ACTUV} which yields estimates for the quasiconformal distortion of sets in terms of the so
 called $h$-contents. This type of contents looks like the usual Hausdorff contents, but they are associated to gauge functions $h$ which depend on the 
point $x$, besides the ``radius'' $r$. They first appeared in \cite{Tolsa-Uriarte}
in the study of the quasiconformal distortion of sets with positive analytic capacity,
and they were further exploited in \cite{ACTUV}, mainly because of their relationship
with the Riesz capacities from non linear potential theory.

Let us remark that previous (weaker) estimates on quasiconformal distortion in terms of Hausdorff measures were obtained by Prause \cite{Prause2}. He proved a result analogous to Theorem \ref{teopri}, with the Hausdorff measure $H^t$ on the left side of \rf{eqhaus} replaced by the Hausdorff measure associated to a gauge function
of the form $r^t\,\ve(r)$, with $\ve(r)$ satisfying $\int_0\ve(r)^\frac{Kt'}{Kt'-t}\,\frac{dr}r<\infty$.

In the paper, as usual, the letter $C$ denotes a constant (often, an absolute constant) that may change 
at different occurrences, while constants with subscript, such as $C_1$, 
retain their values. The notation $A\lesssim B$ means that there is a positive constant $C$ such that $A\leq CB$; and $A\approx B$ means that $A\lesssim B\lesssim A$.

\medskip
\noindent Acknowledgement: I would like to thank Ignacio Uriarte-Tuero for informing me about the problem of proving that $K$-quasiconformal maps 
transform sets of finite $t'$-Hausdorff measure
into sets of finite $t$-Hausdorff measure. 
\bigskip


\section{Preliminaries}

In this section we review some definitions and results from \cite{Tolsa-Uriarte}
and \cite{ACTUV} (mainly dealing with 
the $h$-contents). They will be used to prove Theorem \ref{teopri}.

Let $\mu$ be a finite Borel measure on $\C$. For $0<t< 2$, 
for any ball $B\subset \C$ with radius $r(B)$,
we denote
$\theta_\mu^t(B):=\mu(B)/r(B)^t$. 
Given a parameter $a>0$, we consider the function 
\begin{equation}\label{eqpsia}
\psi_{a,t}(x) = \frac1{|x|^{t+a}+1},\qquad x\in\C.
\end{equation}
For the ball $B=B(x,r)$ we define
\begin{equation}\label{defvex}
\ve_{\mu,a,t}(x,r) = \ve_{\mu,a,t}(B) := \frac1{r^t} \int \psi_{a,t}\Bigl(\frac{y-x}r\Bigr)d\mu(y),
\end{equation}
and we consider the gauge function
\begin{equation}\label{defhx}
h_{\mu,a,t}(x,r) = h_{\mu,a,t}(B) := r^t\ve_{\mu,a,t}(B).
\end{equation}
Notice that $\ve_{\mu,a,t}(B)$ and $h_{\mu,a,t}(B)$ 
can be considered as smooth versions of 
$\theta_\mu^t(B)$ and $\mu(B)$, respectively. One of the advantages of 
$\ve_{\mu,a,t}(x,r)$ over $\theta_\mu(x,r)$ (where, of course, $\theta_\mu^t(x,r):=\theta_\mu^t(B(x,r))$) is that 
 $\ve_{\mu,a,t}(x,2r)\leq C\ve_{\mu,a,t}(x,r)$ for any $x$ and $r>0$, which fails in general for $\theta_\mu^t(x,r)$. Analogously, we have $h_{\mu,a,t}(x,2r)\leq C\,h_{\mu,
 a,t}(x,r)$, while
 $\mu(B(x,r))$ and $\mu(B(x,2r))$ may be very different.

Let $\cB$ denote the family of all closed balls contained in $\C$.
We consider a function $\ve:\cB:\to[0,\infty)$ (for instance, we can take $\ve=\ve_{\mu,a,t}$), and
we define $h(x,r)=\ve(x,r)\,r^t$. We assume that $\ve,h$ are such that $h(x,r)\to0$ as $r\to0$, for all $x\in\C$.
We introduce the measure $H^h$ following Carathéodory's construction (see \cite{mattila}, p.54):
given $0<\delta\leq\infty$ and a set $F\subset\C$, we consider
$$H^h_\delta(F) = \inf\sum_i h(B_i),$$
where the infimum is taken over all coverings $F\subset \bigcup_i B_i$ with balls $B_i$ with radii smaller that $\delta$. Finally, we
define 
$$H^h(F) = \lim_{\delta\to0} H^h_\delta(F).$$
Recall that $H^h$ is a Borel regular measure (see \cite[p.55]{mattila}), although it is not a ``true'' Hausdorff measure.
For the $h$-content, we use the notation $M^h(E):=H_\infty^h(E)$.

We say that the function $\ve(\cdot)$ belongs to $\cG_1$ if it verifies the following properties for all balls 
$B(x,r)$, $B(y,s)$:
there exists a constant $C_2$ such that if $|x-y|\leq 2r$ and $r/2\leq s\leq 2r$, then
\begin{equation}\label{eqeq1}
C_2^{-1}\,\ve(x,r)\leq \ve(y,s)\leq C_2\,\ve(x,r).
\end{equation}
It is easy to check that the function $\ve_{\mu,a,t}$ introduced above belongs to $\cG_1$, for all $0<a<1$. Moreover, 

\begin{lemma} 
For any Borel set $A$, 
$$M^{h_{\mu,a,t}}(A) \geq C(a,t)^{-1}\mu(A).$$
\end{lemma}

See \cite[Lemma 2.2]{Tolsa-Uriarte} for the (easy) proof.

As explained in \cite[Lemma 2.1]{Tolsa-Uriarte}, Frostman Lemma holds for the contents $M^h$: 

\begin{lemma}
Let $\ve\in\cG_1$ and $h(x,r)=r^t\ve(x,r)$. Given a compact set $F\subset \C$, the following holds: $M^h(F)>0$ if and only if
there exists a Borel measure $\nu$ supported on $F$ such that $\nu(B)\leq h(B)$ for any ball $B$. Moreover, one can find $\nu$ such that $\nu(F)\geq c^{-1} M^h(F)$.
\end{lemma}

Given an arbitrary bounded set $A\subset\C$, let $B$ a ball with minimal diameter that contains
$A$. We define $\ve(A):=\ve(B)$. If $B$ is not unique, it does not matter. In this case, for definiteness we choose
the infimum of the values $\ve(B)$ over all balls $B$ with minimal diameter containing $A$, for instance.
Analogously, if $h(x,r)=r^t\,\ve(x,r)$, we define $h(A)$ as the infimum the $h(B)$'s.

Recall that a quasiconformal map $\phi:\C\to\C$ is called principal if $\phi$ is conformal
at $\infty$ (i.e. $\bar\partial \phi$ has compact support) and moreover
$\phi(z) = z (1+ O(z^{-1}))$ as $z\to\infty$.

One of the main technical results from \cite{ACTUV} is the following.

\begin{lemma}\label{mainlem}
 Let $0<t<2$ and $a>0$. Let $\phi:\C\to\C$ be a principal $K$-quasiconformal mapping,
conformal on $\C\setminus \bar \D$.  
 Let $\mu$ be a finite Borel measure on $\C$ and $E\subset B(0,1/2)$.
Denote 
$$
\ve_\phi(x,r):= \ve_{\mu,a,t}(\phi(B(x,r)))^{\frac{t'K}t},\qquad
h_\phi(x,r):= r^{t'}\ve_\phi(x,r).
$$
If $a$ has been chosen small enough (depending only on $t$ and $K$), then we have
$$M^{h_{\mu,a,t}}(\phi(E)) \leq C(K,t) \,M^{h_\phi}(E)^{\frac{t}{t'K}}.$$
\end{lemma}

Let us remark that the case $t=1$ of the preceding result had been obtained previously 
in \cite{Tolsa-Uriarte}.


\bigskip

\section{The proof of Theorem \ref{teopri}}

\subsection{A technical lemma}

\begin{lemma}\label{lemdens}
Let $\mu$ be a finite Borel measure and let $x\in\C$ and $\theta_1>0$ be such that
\begin{equation}\label{eq21}
\frac{\mu(B(x,r))}{r^t} \leq \theta_1 \quad\mbox{ if $0< r\leq \delta$.}
\end{equation}
Then there exists $\delta'>0$ depending only on $\delta,a,t, \theta_1$ and $\mu(\C)$ such that
$$\ve_{\mu,a,t}(x,r)\leq \theta_2 \quad\mbox{ if $0< r\leq \delta'$,}
$$
with $\theta_2$ depending only on $a,t, \theta_1$.
\end{lemma}

\begin{proof}
By the definition of $\ve_{\mu,a,t}$ and $\psi_{a,t}$,
\begin{align*}
\ve_{\mu,a,t}(x,r) &\leq \sum_{j\geq0} \frac1{r^t} \int_{|x-y|\leq 2^jr} \psi_{a,t}
\Bigl(\frac{y-x}r\Bigr)d\mu(y)
 \leq C \sum_{j\geq0} \frac{\mu(B(x,2^jr))}{2^{j(t+a)}\,r^t}.
\end{align*}
 If $2^jr\leq\delta$, we use the estimate
\rf{eq21}. Otherwise, we take into account that
$$\frac{\mu(B(x,2^jr))}{2^{jt}\,r^t} \leq \frac{\mu(\C)}{\delta^t}.$$
So if $N$ denotes the biggest integer such that $2^Nr\leq\delta$, then
\begin{align*}
\ve_{\mu,a,t}(x,r) & \leq C \theta_1 \sum_{0\leq j\leq N} 2^{-ja}+ C 
\frac{\mu(\C)}{\delta^t} \sum_{j\geq N+1} 2^{-ja} \\
& \leq C \biggl(\theta_1 
+ 2^{-Na} \,\frac{\mu(\C)}{\delta^t}\biggr)
\leq C \biggl(\theta_1 
+ \frac{r^a\mu(\C)}{\delta^{t+a}}\biggr).
\end{align*}
If we take $\delta'$ small enough so that 
$$\frac{(\delta')^a\mu(\C)}{\delta^{t+a}} \leq \theta_1,$$
the lemma follows. 
\end{proof}

\bigskip
\subsection{The main lemma}

\begin{lemma}
Let $0<t<2$ and denote 
 $t' = \frac{2t}{2K-Kt+t}$. Let $\phi:\C\to\C$ be 
$K$-quasiconformal.
Assume that $\phi$ is principal and  
conformal outside the unit disk, and that
$E\subset B(0,1/2)$. Then,
 \begin{equation}\label{eqhaus2}
H^t(E) \leq C(K) \,H^{t'}(\phi(E))^{\frac{t}{t'K}}.
\end{equation}
\end{lemma}

\begin{proof}
All the constants $C,C_i$ in this proof are allowed to depend on $t$ and $K$.

To prove \rf{eqhaus2}, we may assume that $H^{t}(E)$ is finite.
Because of the estimates on the upper density of Hausdorff measures (see \cite[p.89]{mattila}), there exists  $\delta>0$ and $F\subset E$ compact such that
$H^{t}(F)\geq H^{t}(E)/2$, and
\begin{equation}\label{eqdens}
\frac{H^{t}(B(x,r)\cap E)}{r^t} \leq 5\quad \mbox{for all $x\in F$ and $0<r\leq \delta.$}
\end{equation}
 Consider the measure $\mu=H^{t}_{|F}$ and the associated $h$-content $M^{h_{\mu,a,t}}$. Recall that we have
\begin{equation}\label{eq1}
M^{h_{\mu,a,t}}(F)\geq C^{-1} \mu(F)= H^{t}(F),
\end{equation}
assuming the parameter $a>0$ above small enough.
By Lemma \ref{mainlem},
\begin{equation}\label{eq2}
M^{h_{\mu,a,t}}(F) \leq C(K,t) \,M^{h_\phi}(\phi(F))^{\frac{t}{t'K}},
\end{equation}
where $M^{h_\phi}$ is the content associated to the gauge function 
$$h_\phi(x,r):= r^{t'}\ve_\phi(x,r),\qquad
\ve_\phi(x,r):= \ve_{\mu,a,t}(\phi^{-1}(B(x,r)))^{\frac{t'K}t}.
$$
By Frostman's lemma we deduce that there exists some measure $\nu$ supported on $\phi(F)$ such that
$\nu(\phi(F))\approx M^{h_\phi}(\phi(F))$ and
\begin{equation}\label{eq61}
\nu(B(x,r))\leq h_\phi(x,r) = r^{t'}\ve_{\mu,a,t}(\phi^{-1}(B(x,r)))^{\frac{t'K}t}.
\end{equation}
From \rf{eqdens} and Lemma \ref{lemdens} we infer that $\ve_{\mu.a.t}(y,s)\leq C_1$ for all $y\in F$ and $0<s<\delta'$, with 
$\delta'$ sufficiently small (depending on $\delta, a,t,\mu(F)$). 
The contant $C_1$ only depends on $a,t$.

By the H\"older continuity of $\phi$, we have
$\ve_{\mu,a,t}(\phi^{-1}(B(x,r))\leq C_1$ for all $x\in \phi(F)$ and $0<r<\delta''$, with $\delta''$ sufficiently small. Indeed, by the definition of $\ve_\phi$, 
there is a ball $B(y,s)\supset
\phi^{-1}(B(x,r))$ with $s\approx\diam\phi^{-1}(B(x,r))$ such that
\begin{equation}\label{id45}
\ve_{\mu,a,t}(\phi^{-1}(B(x,r))= \ve_{\mu,a,t}(B(y,s)).
\end{equation}
Using the H\"older continuity of $\phi$ (see \cite[p.82]{AIM}, for instance), if $B(x,r)\subset \D$, we infer that
$$\frac{s}{\diam(\phi(\D))}\leq C(K)\, \frac{r^{1/K}}{\diam(\D)^{1/K}}.$$ 
Since $\phi$ is principal and conformal outside the unit disk,
$\diam(\phi(\D))\approx\diam(\D)\approx 1$, and thus $s\leq C_2(K)r^{1/K}$.
As a consequence, if $r\leq C_2(K)^{-K}(\delta')^K=:\delta''$, then $s\leq\delta'$, 
and by \rf{id45},
$$\ve_{\mu,a,t}(\phi^{-1}(B(x,r)) \leq C_1,$$
as claimed.

From the preceding estimate and \rf{eq61}, we obtain
$$\nu(B(x,r))\leq C_3r^{t'}\quad\mbox{all $x\in \phi(F)$, $0<r<\delta''$,}$$
with $C_3=C_1^{\frac{t'K}t}$.
This implies that
$H^{t'}(\phi(F))\gtrsim \nu(\phi(F))$. Indeed, if 
$\phi(F)\subset \bigcup_i A_i$, with $\diam(A_i)\leq d\leq \delta''/2$ and $A_i\cap 
\phi(F)\neq\varnothing$,
let $B_i$ be a ball centered on $\phi(F)\cap A_i$ with radius $r(B_i)\leq \diam(A_i)\leq \delta''$.
Then,
$$\sum_i\diam(A_i)^{t'}\geq \sum_i r(B_i)^{t'}\gtrsim \sum_i \nu(B_i) \geq \nu(\phi(F)),$$
and so $H^{t'}_{d}(\phi(F))\gtrsim \nu(\phi(F))$ for all $0<d<\delta''/2$.

From \rf{eq1} and \rf{eq2} it follows that
$$H^t(F) \lesssim M^{h_\phi}(\phi(F))^{\frac{t}{t'K}} \approx 
\nu(\phi(F))^{\frac{t}{t'K}} \lesssim H^{t'}(\phi(F))^{\frac{t}{t'K}}.$$
Since $H^t(E)\approx H^t(F)$ and $H^{t'}(\phi(F))\leq H^{t'}(\phi(E))$, the lemma
follows.
\end{proof}

\bigskip

\subsection{Proof of Theorem \ref{teopri}}

The theorem follows from the preceding lemma by standard arguments in
quasiconformal theory. However,
for completeness we give the details.

We factorize $\phi=\phi_2\circ \phi_1$, where $\phi_1$, $\phi_2$ are both $K$-quasiconformal
maps, with $\phi_1$ principal and conformal on $\C\setminus 2B$, and
$\phi_2$ is conformal on $\phi_1(2B)$. 
Let $g(z)=dz+ b$ be the function that maps the unit disk to $2B$ (so $d=\diam(B)$). 
The function $h=g^{-1}\circ \phi_1\circ g$ verifies the assumptions of the main lemma, and 
thus
$$H^t(g^{-1}(E)) \leq C(K) \,H^{t'}(h(g^{-1}(E)))^{\frac{t}{t'K}}.$$
Since $g$ is an affine map, we have
$$H^t(g^{-1}(E)) = \frac{H^t(E)}{\diam(B)^t},\qquad
H^{t'}(h(g^{-1}(E))) = H^{t'}(g^{-1}(\phi_1(E))) = 
\frac{H^{t'}(\phi_1(E))}{\diam(B)^{t'}}.$$
Using also that  that $\diam(\phi_1(B))\approx\diam(\phi_1(2B))\approx\diam(2B)$
by quasisymmetry and Koebe's distortion theorem, we get
\begin{equation}\label{eqc4}
\frac{H^t(E)}{\diam(B)^t} \leq C(K) \,\left(\frac{H^{t'}(\phi_1(E))}{\diam(\phi_1(B))^{t'}}\right)^{\frac{t}{t'K}}.
\end{equation}

On the other hand, since $\phi_2$ is conformal on $\phi_1(2B)$, by 
Koebe's distortion theorem and quasisymmetry, for each ball $B_0$ contained in 
$B$, 
$$\frac{\diam(\phi_2(\phi_1(B_0)))}{\diam(\phi_2(\phi_1(2B)))} \approx 
\frac{\diam(\phi_1(B_0))}{\diam(\phi_1(2B))}.$$
From this estimate and quasisymmetry again, it is straightforward to check that
$$\frac{H^{t'}(\phi_1(E))}{\diam(\phi_1(B))^{t'}}
\approx \frac{H^{t'}(\phi(E))}{\diam(\phi(B))^{t'}},$$
with constants depending on $K$,
which together with \rf{eqc4} yields \rf{eqhaus}.
\fiproof

\bibliographystyle{alpha}
\bibliography{./refer}

\end{document}